\documentclass[12pt]{article}%
\usepackage{graphicx}
\usepackage[intlimits]{amsmath}
\usepackage{latexsym}
\usepackage{amsfonts}
\usepackage{amssymb}%
\makeatletter
\newfont{\footsc}{cmcsc10 at 8truept}
\newfont{\footbf}{cmbx10 at 8truept}
\newfont{\footrm}{cmr10 at 10truept}
\newtheorem{theorem}{Theorem}

\newtheorem{corollary}[theorem]{Corollary}

\newtheorem{lemma}[theorem]{Lemma}

\newtheorem{problem}[theorem]{Problem}

\newenvironment{proof}[1][Proof]{\noindent{\textbf {#1}  }}  {\hfill$\Box$\bigskip}

\begin{document}

\title{Walks and the spectral radius of graphs}
\author{Vladimir Nikiforov\\Department of Mathematical Sciences, University of Memphis, \\Memphis TN 38152, USA, email: \textit{vnikifrv@memphis.edu}}
\maketitle

\begin{abstract}
Given a graph $G,$ write $\mu\left(  G\right)  $ for the largest eigenvalue of
its adjacency matrix, $\omega\left(  G\right)  $ for its clique number, and
$w_{k}\left(  G\right)  $ for the number of its $k$-walks. We prove that the
inequalities%
\[
\frac{w_{q+r}\left(  G\right)  }{w_{q}\left(  G\right)  }\leq\mu^{r}\left(
G\right)  \leq\frac{\omega\left(  G\right)  -1}{\omega\left(  G\right)  }%
w_{r}\left(  G\right)
\]
hold for all $r>0$ and odd $q>0.$ We also generalize a number of other bounds
on $\mu\left(  G\right)  $ and characterize pseudo-regular and
pseudo-semiregular graphs in spectral terms.



\end{abstract}

\section{Introduction}

Our graph-theoretic notation is standard (e.g., see \cite{Bol98}); in
particular, we assume that graphs are defined on the vertex set $\left\{
1,2,...,n\right\}  =\left[  n\right]  $. Given a graph $G,$ a $k$\emph{-walk}
is a sequence of vertices $v_{1},...,v_{k}$ of $G$ such that $v_{i}$ is
adjacent to $v_{i+1}$ for all $i=1,...,k-1;$ we write $w_{k}\left(  G\right)
$ for the number of $k$-walks in $G$. The eigenvalues of the adjacency matrix
$A\left(  G\right)  $ of $G$ are ordered as $\mu\left(  G\right)  =\mu_{1}%
\geq...\geq\mu_{n}$.

Various bounds of $\mu\left(  G\right)  $ in terms of $w_{k}\left(  G\right)
$ are known; the earliest one, due to Collatz and Sinogowitz \cite{CoSi57},
reads as
\begin{equation}
\mu\left(  G\right)  \geq\frac{2e\left(  G\right)  }{v\left(  G\right)
}=\frac{w_{2}\left(  G\right)  }{w_{1}\left(  G\right)  }. \label{CS}%
\end{equation}
This inequality was strengthened by Hofmeister (\cite{Hof88}, \cite{Hof93}) to%
\begin{equation}
\mu^{2}\left(  G\right)  \geq\frac{1}{v\left(  G\right)  }\sum_{u\in V\left(
G\right)  }d^{2}\left(  u\right)  =\frac{w_{3}\left(  G\right)  }{w_{1}\left(
G\right)  }, \label{Hof}%
\end{equation}
in turn, improved by Yu, Lu, and Tian \cite{YLT04} to
\[
\mu^{2}\left(  G\right)  \geq\frac{w_{5}\left(  G\right)  }{w_{3}\left(
G\right)  },
\]
and by Hong and Zhang \cite{HoZh05} to
\[
\mu^{2}\left(  G\right)  \geq\frac{w_{7}\left(  G\right)  }{w_{5}\left(
G\right)  }.
\]

In this note we prove that, in fact, the inequality
\[
\mu^{r}\left(  G\right)  \geq\frac{w_{q+r}\left(  G\right)  }{w_{q}\left(
G\right)  }%
\]
holds for all $r>0$ and odd $q>0.$

Let $\omega\left(  G\right)  $ be the clique number of $G.$ Wilf \cite{Wil86}
gave the bound
\[
\mu\left(  G\right)  \leq\frac{\omega\left(  G\right)  -1}{\omega\left(
G\right)  }v\left(  G\right)  =\frac{\omega\left(  G\right)  -1}{\omega\left(
G\right)  }w_{1}\left(  G\right)  ,
\]
and Nikiforov \cite{Nik02} showed that
\[
\mu^{2}\left(  G\right)  \leq2\frac{\omega\left(  G\right)  -1}{\omega\left(
G\right)  }e\left(  G\right)  =\frac{\omega\left(  G\right)  -1}{\omega\left(
G\right)  }w_{2}\left(  G\right)  ,
\]
generalizing earlier results in \cite{Cvet72}, \cite{EdEl83}, \cite{Hon95},
\cite{Nos70}, and \cite{Wil86}.

In this note we prove that, in fact, the inequality
\[
\mu^{r}\left(  G\right)  \leq\frac{\omega\left(  G\right)  -1}{\omega\left(
G\right)  }w_{r}\left(  G\right)
\]
holds for every $r\geq1.$

We generalize also a number of other upper and lower bounds on $\mu\left(
G\right)  $ in terms of walks and characterize pseudo-regular
and\ pseudo-semiregular graphs in terms of their eigenvectors.

The rest of the paper is organized as follows. In Section \ref{pre} we recall
some basic notions used further, in Section \ref{lob} we investigate lower
bounds on $\mu\left(  G\right)  ,$ and in Section \ref{UB} we investigate
upper bounds on $\mu\left(  G\right)  .$

\section{\label{pre}Some preliminary results}

Given a graph $G$ and a vertex $u\in V\left(  G\right)  ,$ write
$\Gamma\left(  u\right)  $ for the set of neighbors of $u$ and $w_{k}\left(
u\right)  $ for the number of the $k$-walks starting with $u;$ for every two
vertices $u,v\in V\left(  G\right)  ,$ write $w_{k}\left(  u,v\right)  $ for
the number of the $k$-walks starting with $u$ and ending with $v.$

We state below some basic results related to walks in graphs.

\subsection{\label{NW}The number of $k$-walks in a graph}

Let $G$ be a graph of order $n$ with eigenvalues $\mu_{1}\geq...\geq\mu_{n}$
and $\mathbf{u}_{1},...,\mathbf{u}_{n}$ be corresponding orthogonal unit
eigenvectors. For every $i\in\left[  n\right]  ,$ let $\mathbf{u}_{i}=\left(
u_{i1},...,u_{in}\right)  $ and set $c_{i}=\left(  \sum_{j=1}^{n}%
u_{ij}\right)  ^{2}.$

The number of $k$-walks in $G$ (see, e.g., \cite{CDS80}, p. 44, Theorem 1.10)
is given as follows.

\begin{theorem}
\label{thw}For every $k\geq1,$ $w_{k}\left(  G\right)  =c_{1}\mu_{1}%
^{k-1}+...+c_{n}\mu_{n}^{k-1}.\hfill\square$
\end{theorem}

In particular, for $k=1,$%
\begin{equation}
\sum_{i=1}^{n}c_{i}=n. \label{case1}%
\end{equation}

We also list several equalities that we will use later without reference.%
\[%
\begin{array}
[c]{ll}%
\sum_{u\in V\left(  G\right)  }d^{2}\left(  u\right)  =w_{3}\left(  G\right)
; & \sum_{uv\in E\left(  G\right)  }d\left(  u\right)  d\left(  v\right)
=w_{4}\left(  G\right)  ;\\
\sum_{u\in V\left(  G\right)  }w_{p}^{2}\left(  u\right)  =w_{2p-1}\left(
G\right)  ; & \sum_{u\in V\left(  G\right)  }w_{p}\left(  u\right)
w_{q}\left(  u\right)  =w_{p+q-1}\left(  G\right)  ;\\
\sum_{v\in V\left(  G\right)  }w_{r}\left(  u,v\right)  w_{p}\left(  v\right)
=w_{p+r}\left(  u\right)  ; & \sum_{u,v\in V\left(  G\right)  }w_{r}\left(
u,v\right)  w_{p}\left(  u\right)  w_{q}\left(  v\right)  =w_{p+q+r-2}\left(
G\right)  .
\end{array}
\]

\subsection{The inequality of Motzkin and Straus}

The following result of Motzkin and Straus \cite{MoSt65} will be used in
Section \ref{UB}.

\begin{theorem}
\label{ThMS}For any graph $G$ of order $n$ and real numbers $x_{1},...,x_{n}$
with $x_{i}\geq0,$ $\left(  1\leq i\leq n\right)  ,$ and $x_{1}+...+x_{n}=1,$
\begin{equation}
\sum_{ij\in E\left(  G\right)  }x_{i}x_{j}\leq\frac{\omega\left(  G\right)
-1}{\omega\left(  G\right)  }. \label{MS_in}%
\end{equation}
Equality holds iff the subgraph induced by the vertices corresponding to
nonzero entries of $\mathbf{x}$ is a complete $\omega\left(  G\right)
$-partite graph such that the sum of the $x_{i}$'s in each part is the
same.$\hfill\square$
\end{theorem}

Wilf \cite{Wil86} was the first to apply inequality (\ref{MS_in}) to graph
spectra, obtaining, in particular, the following result.

\begin{theorem}
\label{ThWil}Let $\mathbf{x}=\left(  x_{1},...,x_{n}\right)  $ be an
eigenvector to $\mu\left(  G\right)  $ with $\left\Vert \mathbf{x}\right\Vert
=1.$ Then%
\begin{equation}
\mu\left(  G\right)  =\sum_{ij\in E\left(  G\right)  }x_{i}x_{j}\leq
\frac{\omega\left(  G\right)  -1}{\omega\left(  G\right)  }\left(  \sum
_{i=1}^{n}x_{i}\right)  ^{2}. \label{inWil}%
\end{equation}
$\hfill\square$
\end{theorem}

It is rather entertaining to find the connected graphs for which equality
holds in (\ref{inWil}). We note without a proof that for $G=K_{4n,4n,n}$
equality holds in (\ref{MS_in}) - it is enough to consider the vector
$\mathbf{x}=\left(  x_{1},...,x_{9n}\right)  $ defined as%
\[
x_{i}=\left\{
\begin{array}
[c]{ll}%
\left(  12n\right)  ^{-1/2}, & 1\leq i\leq8n\\
\left(  3n\right)  ^{-1/2}, & 8n<i\leq9n
\end{array}
\right.  .
\]
Here we state only a partial result.

\begin{theorem}
\label{ThWileq}Let $G$ be a connected graph and $\mathbf{x}=\left(
x_{1},...,x_{n}\right)  $ be a unit eigenvector to $\mu\left(  G\right)  $
such that
\[
\mu\left(  G\right)  =\frac{\omega\left(  G\right)  -1}{\omega\left(
G\right)  }\left(  \sum_{i=1}^{n}x_{i}\right)  ^{2}.
\]
Then $G$ is a complete $\omega\left(  G\right)  $-partite graph.
\end{theorem}

\begin{proof}
Since $G$ is connected, $x_{i}>0$ for every $i\in\left[  n\right]  .$ The
assertion follows from the case of equality in (\ref{MS_in}).
\end{proof}

\section{\label{lob}Lower bounds on $\mu\left(  G\right)  $}

Given a graph with no isolated vertices and a vertex $v,$ call the value
$\sum_{v\in\Gamma\left(  u\right)  }d\left(  v\right)  /d\left(  u\right)  $
the \emph{average degree} of $u.$ A graph $G$ with no isolated vertices is called:

- \emph{pseudo-regular} if its vertices have the same average degree;

- \emph{semiregular }if it is bipartite and vertices belonging to the same
part have the same degree;

- \emph{pseudo-semiregular }if it is bipartite and vertices belonging to the
same part have the same average degree.

In this section we first prove Theorem \ref{ThLB} and then show that its
hypothesis cannot be relaxed. Next we describe pseudo-regular and
pseudo-semiregular graphs in terms of their eigenvectors, and finally we
extend two other lower bounds on $\mu\left(  G\right)  .$

The following theorem generalizes results stated in \cite{YLT04} and
\cite{HoZh05}.

\begin{theorem}
\label{ThLB}For every graph $G,$
\begin{equation}
\mu^{r}\left(  G\right)  \geq\frac{w_{q+r}\left(  G\right)  }{w_{q}\left(
G\right)  } \label{min3}%
\end{equation}
for all $r>0$ and odd $q>0.$

If $q>1,$ equality holds in (\ref{min3}) if and only if each component of $G$
has spectral radius $\mu\left(  G\right)  $ and is pseudo-regular or, if $r$
is even, pseudo-semiregular.

If $q=1,$ equality holds in (\ref{min3}) if and only if each component of $G$
has spectral radius $\mu\left(  G\right)  $ and is regular or, if $r$ is even, semiregular.
\end{theorem}

\begin{proof}
Let $v\left(  G\right)  =n.$ Theorem \ref{thw} implies (\ref{min3}) by
\begin{equation}
\frac{w_{q+r}\left(  G\right)  }{w_{q}\left(  G\right)  }=\frac{\sum_{i=1}%
^{n}c_{i}\mu_{i}^{q+r-1}}{\sum_{i=1}^{n}c_{i}\mu_{i}^{q-1}}=\mu^{r}\left(
G\right)  \frac{\sum_{i=1}^{n}c_{i}\left(  \frac{\mu_{i}}{\mu_{1}}\right)
^{q+r-1}}{\sum_{i=1}^{n}c_{i}\left(  \frac{\mu_{i}}{\mu_{1}}\right)  ^{q-1}%
}\leq\mu^{r}\left(  G\right)  . \label{in1}%
\end{equation}

Suppose now that
\begin{equation}
\mu^{r}\left(  G\right)  =\frac{w_{q+r}\left(  G\right)  }{w_{q}\left(
G\right)  }. \label{meq}%
\end{equation}
Assume first that $G$ is connected and let $M$ be the set of all $i\in\left[
2,n\right]  $ such that $c_{i}\neq0$ and $\mu_{i}\neq0.$ We shall show that if
$G$ is nonbipartite, then $M=\varnothing.$ From (\ref{in1}) we find that
\begin{equation}
\sum_{i=2}^{n}c_{i}\left(  \frac{\mu_{i}}{\mu_{1}}\right)  ^{q+r-1}=\sum
_{i=2}^{n}c_{i}\left(  \frac{\mu_{i}}{\mu_{1}}\right)  ^{q-1}, \label{meq1}%
\end{equation}
and so, $\left\vert \mu_{i}\right\vert =\mu_{1}$ for every $i\in M,$
contradicting that $G$ connected and nonbipartite. Hence, $w_{k}\left(
G\right)  =c_{1}\mu_{1}^{k-1}$ for every $k>0.$ In particular, $w_{4}\left(
G\right)  =\sqrt{w_{3}\left(  G\right)  w_{5}\left(  G\right)  },$ and so%
\[
\sum_{u\in V\left(  G\right)  }d\left(  u\right)  w_{3}\left(  u\right)
=w_{4}\left(  G\right)  =\sqrt{w_{3}\left(  G\right)  w_{5}\left(  G\right)
}=\sqrt{\sum_{u\in V\left(  G\right)  }d^{2}\left(  u\right)  \sum_{u\in
V\left(  G\right)  }w_{3}^{2}\left(  u\right)  }.
\]
The condition for equality in Cauchy-Schwarz inequality implies that
$w_{3}\left(  u\right)  /d\left(  u\right)  $ is constant for all vertices
$u,$ i.e., that $G$ is pseudo-regular.

If $q=1,$ then (\ref{meq1}) implies that $c_{i}=0$ for all $1<i\leq n;$ hence
$c_{1}=n,$ $\mu_{1}=w_{2}\left(  G\right)  /n,$ and so $G$ is regular.

Let now $G$ be bipartite. Since the spectrum of $G$ is symmetric with respect
to $0,$ from (\ref{meq1}) it follows that either $M=\varnothing$ or
$M=\left\{  n\right\}  .$ If $M=\varnothing$ (i.e., $c_{n}=0$), the case
reduces to the previous one. If $c_{n}>0$, equality (\ref{meq1}) may hold only
if $r$ is even. Also, we have
\begin{align*}
\sum_{u\in V\left(  G\right)  }d^{2}\left(  u\right)   &  =w_{3}\left(
G\right)  =c_{1}\mu_{1}^{2}+c_{n}\mu_{1}^{2}=\left(  c_{1}+c_{n}\right)
\mu_{1}^{2},\\
\sum_{u\in V\left(  G\right)  }d\left(  u\right)  w_{4}\left(  u\right)   &
=w_{5}\left(  G\right)  =c_{1}\mu_{1}^{4}+c_{n}\mu_{1}^{4}=\left(  c_{1}%
+c_{n}\right)  \mu_{1}^{4},\\
\sum_{u\in V\left(  G\right)  }w_{4}^{2}\left(  u\right)   &  =w_{7}\left(
G\right)  =c_{1}\mu_{1}^{6}+c_{n}\mu_{1}^{6}=\left(  c_{1}+c_{n}\right)
\mu_{1}^{6}.
\end{align*}
Therefore,%
\[
\sum_{u\in V\left(  G\right)  }d^{2}\left(  u\right)  \sum_{u\in V\left(
G\right)  }w_{4}^{2}\left(  u\right)  =\sum_{u\in V\left(  G\right)  }d\left(
u\right)  w_{4}\left(  u\right)  ;
\]
the condition of equality in Cauchy-Schwarz's inequality implies that
$w_{4}\left(  u\right)  /d\left(  u\right)  =w_{4}\left(  v\right)  /d\left(
v\right)  $ for every $u,v\in V\left(  G\right)  .$ We borrow the following
argument from \cite{FMS93}. Letting $u$ to be a vertex of minimum average
degree $w_{3}\left(  u\right)  /d\left(  u\right)  =\delta$ and $v$ be a
vertex of maximum average degree $w_{3}\left(  v\right)  /d\left(  v\right)
=\Delta,$ we see that
\begin{align*}
\frac{w_{4}\left(  u\right)  }{d\left(  u\right)  } &  =\frac{\sum_{t\in
\Gamma\left(  u\right)  }w_{3}\left(  t\right)  }{d\left(  u\right)  }%
\leq\frac{\Delta\sum_{t\in\Gamma\left(  u\right)  }d\left(  t\right)
}{d\left(  u\right)  }=\Delta\delta\\
&  =\frac{\delta\sum_{t\in\Gamma\left(  v\right)  }d\left(  t\right)
}{d\left(  v\right)  }\leq\frac{\sum_{t\in\Gamma\left(  v\right)  }%
w_{3}\left(  t\right)  }{d\left(  v\right)  }\leq\frac{w_{4}\left(  v\right)
}{d\left(  v\right)  }=\frac{w_{4}\left(  u\right)  }{d\left(  u\right)  };
\end{align*}
thus, every vertex of average degree $\delta$ is adjacent only to vertices of
average degree $\Delta$ and vice versa. Since $G$ is connected, it follows
that it is pseudo-semiregular.

If $q=1,$ then
\[
\sum_{i=2}^{n-1}c_{i}\left(  \frac{\mu_{i}}{\mu_{1}}\right)  ^{r}=\sum
_{i=2}^{n-1}c_{i},
\]
and so $c_{i}=0$ for all $1<i<n;$ hence, from (\ref{case1}), $c_{1}+\left(
-1\right)  ^{r}c_{n}=n.$ If $c_{n}=0,$ the case reduces to the previous one.
Otherwise, $r$ is even and so $\mu_{1}^{2}n=w_{3}\left(  G\right)  .$ Since
\[
\mu^{2}\left(  G\right)  >\frac{1}{n}\left\Vert A^{2}\right\Vert _{1}=\frac
{1}{n}w_{3}\left(  G\right)  ,
\]
unless all row sums of $A^{2}\left(  G\right)  $ are equal, we deduce that
$w_{3}\left(  u\right)  $ is constant for every $u,$ and so $G$ is semiregular.

If the graph is not connected, say let $G_{1},...,G_{k}$ be its components, we
have%
\[
\mu^{r}\left(  G\right)  =\frac{\sum_{i=1}^{k}w_{q+r}\left(  G_{i}\right)
}{\sum_{i=1}^{k}w_{q}\left(  G_{i}\right)  }\leq\frac{\sum_{i=1}^{k}\mu
^{r}\left(  G_{i}\right)  w_{q}\left(  G_{i}\right)  }{\sum_{i=1}^{k}%
w_{q}\left(  G_{i}\right)  }\leq\frac{\mu^{r}\left(  G\right)  \sum_{i=1}%
^{k}w_{q}\left(  G_{i}\right)  }{\sum_{i=1}^{k}w_{q}\left(  G_{i}\right)
}\leq\mu^{r}\left(  G\right)  .
\]
Thus, (\ref{meq}) implies that $\mu^{r}\left(  G\right)  =\mu^{r}\left(
G_{i}\right)  =w_{q+r}\left(  G_{i}\right)  /w_{q}\left(  G_{i}\right)  $ for
each component of $G_{i}.$

We omit the straightforward proof of the converse of the case of equality.
\end{proof}

\subsection{The case of even $q$}

Observe that if $G$ is connected and nonbipartite, then the ratio
$w_{q+r}\left(  G\right)  /w_{q}\left(  G\right)  $ tends to $\mu^{r}\left(
G\right)  $ as $q$ tends to infinity. Indeed, from (\ref{in1}) and $\left\vert
\mu_{i}\right\vert /\mu_{1}<1$ holding for every $i=2,...,n,$ we obtain the
following theorem.

\begin{theorem}
For every connected nonbipartite graph $G$ and every $\varepsilon>0,$ there
exists $q_{0}\left(  \varepsilon\right)  $ such that if $q>q_{0}\left(
\varepsilon\right)  $ then
\[
\left(  1-\varepsilon\right)  \frac{w_{q+r}\left(  G\right)  }{w_{q}\left(
G\right)  }\leq\mu^{r}\left(  G\right)  \leq\left(  1+\varepsilon\right)
\frac{w_{q+r}\left(  G\right)  }{w_{q}\left(  G\right)  }%
\]
for every $r>0.$
\end{theorem}

Inequality (\ref{min3}) may fail for $q$ even as shown by the following
example for $q=2k$ and odd $r$. Let $0<a<b$ be integers and $G=K_{a,b}$ be the
complete bipartite graph with parts of size $a$ and $b$. We see that%
\begin{align*}
w_{2k}\left(  G\right)   &  =2a^{k}b^{k},\text{ \ \ \ }\\
w_{2k+r}\left(  G\right)   &  =a\left(  ba\right)  ^{k+\left(  r-1\right)
/2}+b\left(  ba\right)  ^{k+\left(  r-1\right)  /2},\text{\ }\\
\frac{w_{2k+r}\left(  G\right)  }{w_{2k}\left(  G\right)  }  &  =\frac{a+b}%
{2}\left(  ba\right)  ^{\left(  r-1\right)  /2}>\left(  ab\right)  ^{r/2}%
=\mu^{r}\left(  G\right)  .
\end{align*}

Therefore, for bipartite $G,$ $q$ even and $r$ odd, $\mu^{r}\left(  G\right)
$ may differ considerably from $w_{q+r}\left(  G\right)  /w_{q}\left(
G\right)  ,$ no matter how large $q$ is. We are not able to answer the
following natural question.

\begin{problem}
Let $G$ be a connected bipartite graph. Is it true that
\[
\mu^{r}\left(  G\right)  \geq\frac{w_{q+r}\left(  G\right)  }{w_{q}\left(
G\right)  }%
\]
for every even $q\geq2$ and $r\geq2?$
\end{problem}

We also note without a proof that the graph $G=K_{2t,2t,t}$ satisfies $\mu
^{2}\left(  G\right)  <w_{4}\left(  G\right)  /w_{2}\left(  G\right)  .$

\subsection{Characterization of pseudo-regular and pseudo-semiregular graphs}

Write $\mathbf{i}$ for the vector $(1,...,1)\in%
\mathbb{R}
^{n}.$ As a by-product of the proof of Theorem \ref{ThLB} we obtain
characterizations of pseudo-regular and pseudo-semiregular graphs.

\begin{theorem}
If $G$ is a pseudo-regular graph and $\mu_{s}$ is an eigenvalue of $G$ such
that $0<\left\vert \mu_{s}\right\vert <\mu\left(  G\right)  $, then every
eigenvector to $\mu_{s}$ is orthogonal to $\mathbf{i}$. If $G$ has no
bipartite component, then the converse is also true.
\end{theorem}

\begin{proof}
Let $v\left(  G\right)  =n,$ $\mathbf{u}_{1},...,\mathbf{u}_{n}$ be orthogonal
unit eigenvectors of $G$ to $\mu_{1},...,\mu_{n}$ and $c_{1},...,c_{n}$ be as
defined in Section \ref{NW}. Suppose $0<\left\vert \mu_{s}\right\vert <\mu
_{1}$. If $G$ is pseudo-regular, then $w_{4}\left(  G\right)  =\sqrt
{w_{3}\left(  G\right)  w_{5}\left(  G\right)  }$ and so
\[
\sum_{i=1}^{n}c_{i}\left(  \frac{\mu_{i}}{\mu_{1}}\right)  ^{3}=\sqrt
{\sum_{i=1}^{n}c_{i}\left(  \frac{\mu_{i}}{\mu_{1}}\right)  ^{2}\sum_{i=1}%
^{n}c_{i}\left(  \frac{\mu_{i}}{\mu_{1}}\right)  ^{4}}.
\]
The condition for equality in Cauchy-Schwarz's inequality implies that
$\left\vert \mu_{i}/\mu_{1}\right\vert =\mu_{1}/\mu_{1}=1$ whenever
$\left\vert \mu_{i}\right\vert >0$ and $c_{i}>0.$ Hence, $c_{s}=\left(
\sum_{i=1}^{n}u_{si}\right)  ^{2}=0,$ i.e., $\mathbf{u}_{s}$ is orthogonal to
$\mathbf{i}.$

If $G$ has no bipartite component, we see that
\[
w_{k}\left(  G\right)  =\mu^{k-1}\left(  G\right)  \sum_{\mu_{i}=\mu\left(
G\right)  }c_{i}%
\]
for every $k\geq1.$ In particular, $w_{3}\left(  G\right)  \mu^{2}\left(
G\right)  =w_{5}\left(  G\right)  ;$ from the case of equality of Theorem
\ref{ThLB}, we see that $G$ is pseudo-regular, completing the proof.
\end{proof}

\begin{theorem}
\label{thSR}Let $G=G\left(  n\right)  $ be a bipartite graph with eigenvalues
$\mu_{1}\geq...\geq\mu_{n}.$ If $G$ is pseudo-semiregular, then for all
$s\in\left[  n\right]  $ such that $0<\left\vert \mu_{s}\right\vert
<\mu\left(  G\right)  $ every eigenvector to $\mu_{s}$ is orthogonal to
$\mathbf{i}$. If $G$ is connected, the converse is also true.
\end{theorem}

\begin{proof}
Let $\mathbf{u}_{1},...,\mathbf{u}_{n}$ be orthogonal unit eigenvectors to
$\mu_{1},...,\mu_{n},$ and $c_{1},...,c_{n}$ be as defined in Section
\ref{NW}. If $G$ is pseudo-semiregular, $w_{4}\left(  u\right)  /d\left(
u\right)  =w_{4}\left(  v\right)  /d\left(  v\right)  $ for every $u,v\in
V\left(  G\right)  .$ Letting $t=w_{4}\left(  u\right)  /d\left(  u\right)  ,$
we see that $w_{2k+3}\left(  u\right)  =tw_{2k+1}\left(  u\right)  $ for every
integer $k>0.$ Hence, $t=\mu^{2}\left(  G\right)  $ and so $w_{5}\left(
G\right)  =\mu^{2}\left(  G\right)  w_{3}\left(  G\right)  ,$ implying in turn%
\[
\sum_{i=2}^{n-1}c_{i}\left(  \frac{\mu_{i}}{\mu_{1}}\right)  ^{4}=\sum
_{i=2}^{n-1}c_{i}\left(  \frac{\mu_{i}}{\mu_{1}}\right)  ^{2}.
\]
We see that $c_{s}=0$ for every $s$ such that $0<\left\vert \mu_{s}\right\vert
<\mu\left(  G\right)  ,$ and so $\mathbf{u}_{s}$ is orthogonal to $\mathbf{i}$.

If $G$ is connected and for every $s$ such that $0<\left\vert \mu
_{s}\right\vert <\mu\left(  G\right)  $ every eigenvector to $\mu_{s}$ is
orthogonal to $\mathbf{i}$, then
\[
w_{k}\left(  G\right)  =\left(  c_{1}+\left(  -1\right)  ^{k-1}c_{n}\right)
\mu^{k-1}\left(  G\right)  ,
\]
for every integer $k>0.$ In particular, $w_{5}\left(  G\right)  =\mu
^{2}\left(  G\right)  w_{3}\left(  G\right)  ;$ from the case of equality in
Theorem \ref{ThLB}, it follows that $G$ is pseudo-semiregular.
\end{proof}

\subsection{More lower bounds}

A common device for finding lower bounds on $\mu\left(  G\right)  $ is the
Rayleigh principle applied with carefully chosen vectors.

Let $p\geq0,$ $r\geq1$ be integers and $G$ be a graph\ of order $n$ with no
isolated vertices. Setting $x_{i}=w_{p}\left(  i\right)  /\sqrt{w_{2p-1}%
\left(  G\right)  }$ for all $i\in\left[  n\right]  $ and letting
$\mathbf{x}=\left(  x_{1},...,x_{n}\right)  ,$ the Rayleigh principle gives
another proof of inequality (\ref{min3}) by%
\[
\mu^{r}\left(  G\right)  \geq\left\langle A^{r}\left(  G\right)
\mathbf{x},\mathbf{x}\right\rangle =\frac{1}{w_{2p-1}\left(  G\right)  }%
\sum_{u,v\in V\left(  G\right)  }w_{r+1}\left(  u,v\right)  w_{p}\left(
u\right)  w_{p}\left(  v\right)  =\frac{w_{2p+r-1}\left(  G\right)  }%
{w_{2p-1}\left(  G\right)  }.
\]

Set now $y_{i}=\sqrt{w_{p}\left(  i\right)  /w_{p}\left(  G\right)  }$ for all
$i\in\left[  n\right]  $ and let $\mathbf{y}=\left(  y_{1},...,y_{n}\right)
.$ By the Rayleigh principle we obtain the following general bound%
\begin{equation}
\mu^{r}\left(  G\right)  \geq\left\langle A^{r}\left(  G\right)
\mathbf{y},\mathbf{y}\right\rangle =\frac{1}{w_{p}\left(  G\right)  }%
\sum_{u,v\in V\left(  G\right)  }w_{r+1}\left(  u,v\right)  \sqrt{w_{p}\left(
u\right)  w_{p}\left(  v\right)  }. \label{FMS1}%
\end{equation}
Since by Cauchy-Schwarz's inequality we have%
\begin{align*}
\sum_{u,v\in V\left(  G\right)  }w_{r+1}\left(  u,v\right)  \sqrt{w_{p}\left(
u\right)  w_{p}\left(  v\right)  }\sum_{u,v\in V\left(  G\right)  }%
\frac{w_{r+1}\left(  u,v\right)  }{\sqrt{w_{p}\left(  u\right)  w_{p}\left(
v\right)  }}  &  \geq\left(  \sum_{u,v\in V\left(  G\right)  }w_{r+1}\left(
u,v\right)  \right)  ^{2}\\
&  =w_{r+1}^{2}\left(  G\right)  ,
\end{align*}
inequality (\ref{FMS1}) implies also that
\begin{equation}
\mu^{r}\left(  G\right)  \sum_{u,v\in V\left(  G\right)  }\frac{w_{r+1}\left(
u,v\right)  }{\sqrt{w_{p}\left(  u\right)  w_{p}\left(  v\right)  }}\geq
\frac{w_{r+1}^{2}\left(  G\right)  }{w_{p}\left(  G\right)  }. \label{FMS2}%
\end{equation}

Setting $p=2,r=1,$ we obtain the following inequalities proved by Favaron,
Mah\'{e}o, and Sacl\'{e} \cite{FMS93}, and in a wider context also by Hoffman,
Wolfe, and Hofmeister \cite{HWH95},
\begin{align}
\mu\left(  G\right)   &  \geq\frac{1}{2e\left(  G\right)  }\sum_{uv\in
E\left(  G\right)  }\sqrt{d\left(  u\right)  d\left(  v\right)  }%
,\label{Hof1}\\
\mu\left(  G\right)   &  \geq\frac{2e\left(  G\right)  }{\sum_{uv\in E\left(
G\right)  }\frac{1}{\sqrt{d\left(  u\right)  d\left(  v\right)  }}}.
\label{Hof2}%
\end{align}
As shown in \cite{FMS93} and \cite{HWH95} equality holds in (\ref{Hof1}) and
(\ref{Hof2}) iff $G$ is regular or semiregular. The case of equality in
(\ref{FMS1}) and (\ref{FMS2}) is an open question.

\section{\label{UB} Upper bounds on $\mu\left(  G\right)  $}

In this section we present two general upper bounds on $\mu\left(  G\right)
$. Theorem \ref{thUB} below gives the first bound in terms of the clique
number and the number of walks. The bound of the second type is given in
Section \ref{MUP}. The proof of Theorem \ref{thUB} relies on two simple
preliminary results.

\begin{lemma}
\label{leMSw0}For every $r>0$ and every graph $G,$
\[
w_{2r}\left(  G\right)  \leq\frac{\omega\left(  G\right)  -1}{\omega\left(
G\right)  }w_{r}^{2}\left(  G\right)  .
\]

\end{lemma}

\begin{proof}
Indeed, we have%
\[
w_{2r}\left(  G\right)  =\sum_{uv\in E\left(  G\right)  }w_{r}\left(
u\right)  w_{r}\left(  v\right)  \leq\frac{\omega\left(  G\right)  -1}%
{\omega\left(  G\right)  }\left(  \sum_{u\in V\left(  G\right)  }w_{r}\left(
u\right)  \right)  ^{2}=\frac{\omega\left(  G\right)  -1}{\omega\left(
G\right)  }w_{r}^{2}\left(  G\right)  .
\]

\end{proof}

Applying Lemma \ref{leMSw0} several times, we generalize it as follows.

\begin{corollary}
\label{leMSw}For every graph $G$ and $k,r>0,$
\[
\frac{\omega\left(  G\right)  -1}{\omega\left(  G\right)  }w_{2^{k}r}\left(
G\right)  \leq\left(  \frac{\omega\left(  G\right)  -1}{\omega\left(
G\right)  }w_{r}\left(  G\right)  \right)  ^{2^{k}}.
\]

\end{corollary}

We are ready now to prove the main result of this section.

\begin{theorem}
\label{thUB}For every graph $G$ and $r\geq1,$%
\begin{equation}
\mu^{r}\left(  G\right)  \leq\frac{\omega\left(  G\right)  -1}{\omega\left(
G\right)  }w_{r}\left(  G\right)  \label{upb}%
\end{equation}

\end{theorem}

\begin{proof}
Clearly it suffices to prove inequality (\ref{upb}) for connected graphs. We
shall assume first that $G$ is nonbipartite. Assume that (\ref{upb}) fails,
i.e.,
\[
\mu^{r}\left(  G\right)  >\left(  1+c\right)  \frac{\omega\left(  G\right)
-1}{\omega\left(  G\right)  }w_{r}\left(  G\right)
\]
for some $G,$ $r>0,$ $c>0$. Then, by Corollary \ref{leMSw}, for every $k>0,$%
\begin{equation}
\mu^{2^{k}r}\left(  G\right)  >\left(  \left(  1+c\right)  \frac{\omega\left(
G\right)  -1}{\omega\left(  G\right)  }w_{r}\left(  G\right)  \right)
^{2^{k}}\geq\left(  1+c\right)  ^{2^{k}}\frac{\omega\left(  G\right)
-1}{\omega\left(  G\right)  }w_{2^{k}r}\left(  G\right)  . \label{asin1}%
\end{equation}
Note that Theorem \ref{thw} implies that for every $\varepsilon,$
\begin{equation}
c_{1}\mu^{q-1}\left(  G\right)  <\left(  1+\varepsilon\right)  w_{q}\left(
G\right)  \label{upb1}%
\end{equation}
for all sufficiently large $q$. Hence, for $q=2^{k}r$ and $k$ sufficiently
large, Theorem \ref{ThWil} and inequality (\ref{upb1}) imply that
\[
\mu^{2^{k}r}\left(  G\right)  \leq\frac{\omega\left(  G\right)  -1}%
{\omega\left(  G\right)  }c_{1}\mu^{2^{k}r-1}\left(  G\right)  <\left(
1+\varepsilon\right)  \frac{\omega\left(  G\right)  -1}{\omega\left(
G\right)  }w_{2^{k}r}\left(  G\right)  ,
\]
contradicting (\ref{asin1}).

Finally we have to prove (\ref{upb}) for bipartite $G.$ Then $\omega\left(
G\right)  =2,$ so we have to prove that $\mu^{r}\left(  G\right)  \leq
w_{r}\left(  G\right)  /2$ for every $r\geq2.$ If $r$ is odd, Theorem
\ref{ThWil} and Theorem \ref{thw} imply%
\[
\mu^{r}\left(  G\right)  \leq\frac{1}{2}c_{1}\mu_{1}^{r-1}\leq\frac{1}{2}%
w_{r}\left(  G\right)  .
\]

Let now $r$ be even. Write $cw_{k}\left(  G\right)  $ for the number of
\emph{closed walks} on $k$ vertices in $G$ (i.e., $k$-walks with the same
start and end vertex.) It is known that
\begin{equation}
cw_{k+1}\left(  G\right)  =tr\left(  A^{k}\left(  G\right)  \right)  =\mu
_{1}^{k}+...+\mu_{n}^{k}. \label{cwk}%
\end{equation}
The spectrum of bipartite graphs is symmetric with respect to $0,$ thus
$2\mu^{r}\left(  G\right)  \leq cw_{r+1}\left(  G\right)  \leq w_{r}\left(
G\right)  ,$ completing the proof.
\end{proof}

\begin{theorem}
Suppose that $G$ is graph such that equality holds in (\ref{upb}) for some
$r\geq1$. If $r=1$, then $G$ is a regular complete $\omega\left(  G\right)
$-partite graph. If $r>1$, then $G$ has a single nontrivial component $G_{1}.$
If $\omega\left(  G\right)  >2$, then $G_{1}$ is a regular complete
$\omega\left(  G\right)  $-partite graph. If $\omega\left(  G\right)  =2$,
then $G_{1}$ is a complete bipartite graph, and if $r$ is odd, then $G_{1}$ is regular.
\end{theorem}

\begin{proof}
Assume
\begin{equation}
\mu^{r}\left(  G\right)  =\frac{\omega\left(  G\right)  -1}{\omega\left(
G\right)  }w_{r}\left(  G\right)  \label{case2}%
\end{equation}
and let $c_{i}$ be defined as in Section \ref{NW}.

If $r=1$ then
\[
\frac{\omega\left(  G\right)  -1}{\omega\left(  G\right)  }v\left(  G\right)
=\mu\left(  G\right)  \leq\sqrt{2\frac{\omega\left(  G\right)  -1}%
{\omega\left(  G\right)  }e\left(  G\right)  };
\]
from the case of equality in Tur\'{a}n's theorem (see, e.g., \cite{Bol98}) it
follows that $G$ is regular complete $\omega\left(  G\right)  $-partite graph.

Assume now $r\geq2$; let $G_{1}$ be a component of $G$ with $\mu\left(
G\right)  =\mu\left(  G_{1}\right)  $. If $G_{2}$ is another nontrivial
component of $G$, then
\[
\mu^{r}\left(  G_{1}\right)  =\mu^{r}\left(  G\right)  =\frac{\omega\left(
G\right)  -1}{\omega\left(  G\right)  }\left(  w_{r}\left(  G_{2}\right)
+w_{r}\left(  G_{1}\right)  \right)  >\frac{\omega\left(  G\right)  -1}%
{\omega\left(  G\right)  }w_{r}\left(  G_{1}\right)  ,
\]
a contradiction; thus $G_{1}$ is the only nontrivial component of $G.$ We also
see that the equality (\ref{case2}) holds for $G_{1},$ so for simplicity we
shall assume that $G$ is connected. From Corollary \ref{leMSw} and
(\ref{case2}) we deduce that
\begin{equation}
\mu^{2^{k}r}\left(  G\right)  =\frac{\omega\left(  G\right)  -1}{\omega\left(
G\right)  }w_{2^{k}r}\left(  G\right)  =\frac{\omega\left(  G\right)
-1}{\omega\left(  G\right)  }\mu^{2^{k}r-1}\left(  G\right)  \sum_{i=1}%
^{n}c_{i}\left(  \frac{\mu_{i}}{\mu_{1}}\right)  ^{2^{k}r-1} \label{case3}%
\end{equation}
for every integer $k>0.$ Assume $G$ is nonbipartite; therefore, $\left\vert
\mu_{n}\left(  G\right)  \right\vert <\mu\left(  G\right)  $ and, letting $k$
tend to infinity, we find that
\[
\mu\left(  G\right)  =\frac{\omega\left(  G\right)  -1}{\omega\left(
G\right)  }c_{1}.
\]
From Theorem \ref{ThWileq} it follows that $G$ is a complete $\omega\left(
G\right)  $-partite graph, and thus $G$ has no positive eigenvalues other than
$\mu\left(  G\right)  $. Hence, from (\ref{case3}), any $c_{i}$ corresponding
to a negative eigenvalue must be $0.$ Therefore,%
\[
n=w_{1}\left(  G\right)  =c_{1}\mu\left(  G\right)  =\frac{\omega\left(
G\right)  }{\omega\left(  G\right)  -1}\mu\left(  G\right)  ,
\]
a case that is settled above.

Let now $G$ be bipartite. If $r$ is odd, we have%
\[
\mu^{r}\left(  G\right)  =\frac{1}{2}w_{r}\left(  G\right)  =\frac{1}{2}%
\sum_{i=1}^{n}c_{i}\left(  \mu_{i}\right)  ^{r-1};
\]
so, by Theorem \ref{ThWil} $c_{1}=1/2.$ Moreover, either $c_{i}=0$ of $\mu
_{i}=0$ for $i=2,...,n.$ We have again
\[
n=w_{1}\left(  G\right)  =c_{1}\mu\left(  G\right)  =2\mu\left(  G\right)  ,
\]
implying that $G$ is a regular complete bipartite graph.

For even $r$ we have
\[
2\mu^{r}\left(  G\right)  \leq cw_{r+1}\left(  G\right)  \leq w_{r}\left(
G\right)  =2\mu^{r}\left(  G\right)  ,
\]
and, in view of (\ref{cwk}), we conclude that $G$ has only two nonzero
eigenvalues - $\mu_{1}$ and $\mu_{n}.$ Hence, in our case, Smith's theorem
implies that $G$ is a complete bipartite graph.
\end{proof}

\subsection{\label{MUP}More upper bounds}

It is known that the Perron root of a nonnegative matrix does not exceed its
maximal row sum. This idea has been exploited to obtain the following bounds
\begin{align}
\mu\left(  G\right)   &  \leq\max_{u\in V\left(  G\right)  }\sqrt{w_{3}\left(
u\right)  },\label{lin0}\\
\mu\left(  G\right)   &  \leq\max_{u\in V\left(  G\right)  }\frac{w_{3}\left(
u\right)  }{d\left(  u\right)  },\label{lin1}\\
\mu\left(  G\right)   &  \leq\max_{uv\in E\left(  G\right)  }\sqrt{d\left(
u\right)  d\left(  v\right)  },\label{lin2}\\
\mu\left(  G\right)   &  \leq\max_{uv\in E\left(  G\right)  }\sqrt{\frac
{w_{3}\left(  u\right)  w_{3}\left(  v\right)  }{d\left(  u\right)  d\left(
v\right)  }}. \label{lin3}%
\end{align}

Inequalities (\ref{lin0}) and (\ref{lin1}) are proved in \cite{FMS93},
inequality (\ref{lin2}) is proved in \cite{BeZh01}, and inequality
(\ref{lin3}) in \cite{DaKu04}. As an attempt to interrupt this monotonic
sequence we propose the following general result.

\begin{theorem}
For every integers $p\geq1,$ $r\geq1$ and any graph $G,$
\[
\mu^{r}\left(  G\right)  \leq\max_{u\in V\left(  G\right)  }\frac
{w_{r+p}\left(  u\right)  }{w_{p}\left(  u\right)  }.
\]

\end{theorem}

\begin{proof}
Set $b_{ii}=w_{p}\left(  i\right)  $ for each $i\in\left[  n\right]  $ and let
$B$ be the diagonal matrix with main diagonal $\left(  b_{11},...,b_{nn}%
\right)  .$ Since $B^{-1}A^{r}\left(  G\right)  B$ has the same spectrum as
$A^{r}\left(  G\right)  ,$ $\mu^{r}\left(  G\right)  $ is bounded from above
by the maximum row sum of $B^{-1}A^{r}\left(  G\right)  B$ - say the sum of
the $k$th row - and so,
\[
\mu^{r}\left(  G\right)  \leq\sum_{v\in V\left(  G\right)  }w_{r}\left(
k,v\right)  \frac{w_{p}\left(  v\right)  }{w_{p}\left(  k\right)  }%
=\frac{w_{r+p}\left(  k\right)  }{w_{p}\left(  k\right)  }\leq\max_{u\in
V\left(  G\right)  }\frac{w_{r+p}\left(  u\right)  }{w_{p}\left(  u\right)
}.
\]

\end{proof}

Setting $p=1,$ $r=2,$ we obtain (\ref{lin0}); the case $p=2,$ $r=1$ implies
(\ref{lin1}). Furthermore, (\ref{lin2}) follows from (\ref{lin0}) by%
\[
\mu^{2}\left(  G\right)  \leq\max_{u\in V\left(  G\right)  }w_{3}\left(
u\right)  =\max_{uv\in E\left(  G\right)  }d\left(  u\right)  \left(  \frac
{1}{d\left(  u\right)  }\sum_{v\in\Gamma\left(  u\right)  }d\left(  v\right)
\right)  \leq\max_{uv\in E\left(  G\right)  }d\left(  u\right)  d\left(
v\right)  ,
\]
and (\ref{lin3}) follows by
\begin{align*}
\mu^{2}\left(  G\right)   &  \leq\max_{u\in V\left(  G\right)  }\frac
{w_{4}\left(  u\right)  }{d\left(  u\right)  }=\max_{u\in V\left(  G\right)
}\frac{w_{3}\left(  u\right)  }{d\left(  u\right)  }\frac{w_{4}\left(
u\right)  }{w_{3}\left(  u\right)  }\leq\max_{u\in V\left(  G\right)  }%
\frac{w_{3}\left(  u\right)  }{d\left(  u\right)  }\left(  \frac{\sum
_{v\in\Gamma\left(  u\right)  }w_{3}\left(  v\right)  }{\sum_{v\in
\Gamma\left(  u\right)  }d\left(  v\right)  }\right) \\
&  \leq\max_{u\in V\left(  G\right)  }\frac{w_{3}\left(  u\right)  }{d\left(
u\right)  }\left(  \frac{1}{d\left(  u\right)  }\sum_{v\in\Gamma\left(
u\right)  }\frac{w_{3}\left(  v\right)  }{d\left(  v\right)  }\right)
\leq\max_{uv\in E\left(  G\right)  }\frac{w_{3}\left(  u\right)  }{d\left(
u\right)  }\frac{w_{3}\left(  v\right)  }{d\left(  v\right)  }%
\end{align*}
with plenty of room.\bigskip

\textbf{Acknowledgement }I am grateful to the referee for his valuable remarks.

\end{document}